\documentclass[3p]{elsarticle}

\usepackage{hyperref}

\usepackage[utf8]{inputenc}
\usepackage{tikz}
\usepackage{amsmath,amssymb,amsthm}
\usepackage{bbm}

\theoremstyle{plain}
\newtheorem{thm}{Theorem}[section]

\theoremstyle{definition}

\theoremstyle{remark}
\newtheorem*{rem}{Remark}

\foreach \a in {q,w,e,r,t,y,u,i,o,p,a,s,d,f,g,h,j,k,l,z,x,c,v,b,n,m,Q,W,E,R,T,Y,U,I,O,P,A,S,D,F,G,H,J,K,L,Z,X,C,V,B,N,M}{
  \expandafter\xdef\csname v\a\endcsname {
		{\noexpand\boldsymbol{\a}}
	}
}

\newcommand{\vone}{{\boldsymbol{1}}}

\newcommand{\ind}[1]{{\mathbbm{1}}_{#1}}

\allowdisplaybreaks[4]

\begin{document}

\begin{frontmatter}

\title{Opinion Dynamics on Graphon: the piecewise constant case}

\author[esp,adamss]{Giacomo Aletti}
\ead{giacomo.aletti@unimi.it}
\author[esp,adamss,cnr]{Giovanni Naldi\corref{cor1} \fnref{fn2}}
\ead{giovanni.naldi@unimi.it}

\cortext[cor1]{Corresponding author}

\address[esp]{Department of Environmental Science and Policy, Universit\`{a} degli Studi di Milano, \\ via Celoria 2, 20133 Milano, Italy}
\address[adamss]{ADAMSS Center, Universit\`{a} degli Studi di Milano, via Celoria 2, 20133 Milano, Italy}
\address[cnr]{IMATI-CNR, Via Ferrata, 5/a, 27100 Pavia, Italy}





\begin{abstract}
The study of network supported opinion
dynamics in large groups of autonomous agents is attracting an increasing interest during the last years.
In this paper,  
we proposed the use of the recent graphon theory 
to model and simulate an interacting system
describing the evolution of individual opinions over an arbitrary size networks. 
Specifically, we prove the existence and uniqueness of the limit problem that approximates
a very large networks made by homogeneous groups of agents.
The significant new example is the mean field analysis deduced from the graphon limit systems 
in the case of piecewise constant graphon.
\end{abstract}

\end{frontmatter}


\section{Introduction}
Opinion dynamics is an active research area that deals with the evolution of opinions through the social interaction between a group of
individuals or agents.  In order to understand this dynamics we have to define the individual's opinion, we
should explain how the agents interact, and we should include any external influence factor \cite{IJMOC2020,PT2013}.

Opinions are described by both discrete or continuous quantities. Some examples about the first models are Sznajd model \cite{sznajd2000opinion}, voter model \cite{clifford1973model}, majority rule model \cite{galam1986}, and the Latan\`{e} model of the social impact.
In continuous models, the distribution of opinions are used to be represented by real numbers.  
Examples from this class of models include the DeGroot model \cite{Degroot1974118},
the FJ model \cite{Friedkin1990193}, the Hegselmann-Krause model \cite{Hegselmann2002}, and the
Deffuant-Weisbuch model \cite{Deffuant2000}.

Opinion formation is a complex process affected by the interplay of different elements, including the individual predisposition, the influence of positive and negative
interactions with other individual, the information each individual is exposed to, and many others. A social network can be considered as a set of people (or even groups of people) that participate and interact sharing different kinds of information with the purpose of friendship, decision making,  business exchange or marketing. Then,  a crucial aspect is
the analysis of the different structures in the network to understand what may either facilitate or not the formation of collective belief.
A natural way to represent these interactions is a network which can be described by a weighted graph $G^M=(V_M,E_M,B_M)$, $M>1$, where the nodes' set $V_M=\{1,2,...,M\}$ represents the  set of individuals (agents),  and an edge set $E_M\subseteq V_M \times V_M$, representing pairwise interactions between each pair of agents \cite{MYE2019,Noorazar2020,Benf6,Benf7}. The edges in the set $E_M$  are unoriented, $e=(i,j)=(j,i),\,i,j\in V$. 
In this paper we assign a weight $B_M(e)$ to each edge $e\in E$, where $B_M:E \rightarrow [-K,K]$, with $K>0$. 
$B_M((i,j))$ models the persuasiveness of the agent $i$ with respect to the agent $j$. 
We point out that we allow $B_M$ to be negative, that means that the one agent might
oppose to another one ending in a farther opinions, together with agents that mediate their
opinion ($B_M$ positive). This model may be seen as reacher than pure spreading or concentration models (see , e.g., \cite{ANT07}).
External factors are
made by media or by the so-called opinion leaders, who are more active agents that transmit information without being affected by other agents. 

An equivalent representation of the graph $G^M$ is obtained by the weighted adjacency matrix $B^M$ which is defined to be the square matrix $M\times M$ such that an element $B^M_{ij}$ is $B_M(e)$ when there is an edge $e$ from node $i$ to node $j$, and zero otherwise. To each node $i$ we associate an opinion value $u_i^M(t)\in \mathbb{R}$ at continuous time $t\geq 0$.  Then, we model the opinion dynamics by the following first order system of differential equations (see e.g. \cite{Hegselmann2002})

\begin{equation}\label{syst:op:dyn}
\dot{u}^M_i(t)=\frac{1}{M} \sum_{j=1}^{M} B^M_{ij} \left( u_j^M(t) - u_i^M(t) \right),\,\,i=1,2,...,M,
\end{equation}
with an initial opinion value $u_i^M(0)=u_{i0}^M$, and where we normalized the weights with the factor $1/M$.

In order to study the dynamics of \eqref{syst:op:dyn} for a very large network, as in the case of the recent social media with virtual interactions between people and services,  we will use the so-called graphon theory which is a new paradigm for understanding the continuum limit of a graph sequences when the number of the nodes $M$ goes to infinity  \cite{Borgs2008,Borgs2012}. 
To consider a graphon associated to the graph $G^M$ we partition the interval $[0,1]$ into $M$ isometric intervals $I_j^M := ((j-1)/M,j/M]$, $j=1,2,...,M$. The graphon $W_{G^M}(x,y):[0,1]\times [0,1]\rightarrow\mathbb{R}$ corresponding to the graph $G^M$ is defined by setting $W_{G^M}(x,y)=B^M_{ij}$ if $(x,y)\in I_i^M \times I_j^M$ (this also called pixel diagram,
see \cite{Gao20}). In this work  a graphon is a bounded symmetric Lebesgue measurable functions $W:[0,1]\times [0,1]\rightarrow [-K,K]$, $K>0$, which can be interpreted as weighted graphs on the node set $[0,1]$.
In the space of graphons we define the so-called cut-norm

\[
\| W \|_\Box = \sup_{S,T\subseteq [0,1]}\left| \int_{S\times T} W(x,y) dx dy  \right| ,
\]
with the supremum taken over all measurable subsets $S$ and $T$ of $[0,1]$. The following inequalities hold between norms on a graphon $W$:
$\| W \|_\Box\leq \| W \|_1\leq \| W \|_2\leq \| W \|_\infty \leq const.$ Now the vector function of the opinion can be identify with the piecewise constant function
$u_M(x,t)=\sum_{i=1}^{M} u_j(t) \ind{I_j^M}(x)$, where $\ind{I_j^M}$ is the characteristic function of the interval $I_j^M$. We assume that 
the entries of the $M$-th adjacency matrix satisfy

\[
\sup_{M\in\mathbb{N}}\,\sup_{j=1,2,...,M}\, \frac{1}{M} \sum_{k=1}^{M} B^M_{ij} \leq C,
\]
for a suitable constant $C>0$, and the graphons $W_{G^M}$ converge to a graphon $W$ in the cut-norm (see \cite[Section III]{Gao20} and the references
therein). 
Then it is possible to show (see, e.g., \cite[Theorem~4]{Gao20}) that the 
sequence of opinion vector functions $u_M(x,t)$ converges in $L^2([0,T],L^2([0,1]))$ to the solution of the following  graphon Chauchy problem

\begin{equation}\label{eq:GraphonChaucy}
\begin{cases}
\frac{ \partial  u}{\partial t}(x,t) = \int_{[0,1]} W(x,y) (u(y,t) - u(x,t)) dy
\\
u(x,0) = u_0(x)
\end{cases}
\end{equation}
In this case we can understand the opinion dynamics on a very large graph $G_M$ by using the analytical properties of the solution of the  mean-field equation \eqref{eq:GraphonChaucy}. Therefore with the graphon approach we can simplify the study of the behavior of the solution of \eqref{syst:op:dyn} 
when a simpler analytical solution may be achieved by \eqref{eq:GraphonChaucy}.  In the following we will consider as a significant example the case of the step function graphons.
This approach might be used in the clustering problems for image analysis described in \cite{Benf9,Benf11} or in
dynamical systems on networks as in \cite{Benf10,KVM2018}.

\section{Theoretical analysis of the step function graphons}
We consider the case of $N$ subgroups (parties, homogeneous subset of people), and we suppose different communication between people of each subgroup with respect to people of different groups. A good approximation of this realistic situation can be obtained by considering step function graphons.
The dynamics of the system leads to the study of the possible reinforcement of groups or to the fragmentation of opinion.
Let $N$ be fixed; we define a \emph{partition} $\vi$ of length $N$ as
$\vi = (i_0, , i_1 , i_2 , \ldots , i_{N-1} , i_N) $, where
$0 = i_0 < i_1 < i_2 < \cdots < i_{N-1} \leq i_N = 1 $.
For any interval $I \subseteq [0,1]$, we call $\mathcal{P}_{I}$ the set of real functions on $[0,1]$ constant in $I$ and $0$ outside:
$\mathcal{P}_{I} = \{ u = a\ind{I}, a \in \mathbb{R} \}$.
Given a partition $\vi$, let $\mathcal{P}_{\vi}$ be the set of real functions on $[0,1]$ constant on $\vi$:
 $\mathcal{P}_{\vi} = \{ u = \sum_{j=1}^N a_j \ind{( i_{j_1-1},i_{j_1}]} , a_j \in \mathbb{R} \}$.
Moreover, we call $\mathcal{P}_{\vi\times\vi} $ the set of real functions on $[0,1]\times [0,1]$
constant on the rectangular subdivisions, and we assume that $W \in \mathcal{P}_{\vi\times\vi}$: for $(x,y) \in [0,1]^2$

\[
W(x,y) = \sum_{j_1,j_2=1}^N b_{j_1j_2} \ind{( i_{j_1-1},i_{j_1}]}(x) \ind{( i_{j_2-1},i_{j_2}]}(y).
\]
From now on, we assume that the graphon kernel $W \in \mathcal{P}_{\vi\times\vi} $.
For any interval $I \subseteq [0,1]$, we call $L^2_0 (I)$ the set of real functions $u\in L^2([0,1])$ which are zero outside $I$ and
with zero mean:
$L^2_0 (I) = \{u\in L^2([0,1]) \colon u = u \ind{I}, \int_{[0,1]} u = \int_{I} u = 0 \}$.
It is obvious that
$L^2(0,1) = \mathcal{P}_{\vi} \oplus_{j=1}^N L^2_0 ((i_{j-1},i_{j}]) $, where the sum involves orthogonal closed subspaces of $L^2(0,1)$.
In fact, for any $u \in L^2([0,1])$, we have

\begin{multline}\label{eq:superposit}
u (x)= \sum_{i=1}^N \ind{(i_{j-1},i_{j}]} (x)\Big( u(x) - \frac{\int_{(i_{j-1},i_{j}]} u(y)dy}{i_{j}-i_{j-1}}
+ \frac{\int_{(i_{j-1},i_{j}]} u(y)dy}{i_{j}-i_{j-1}}  \Big)
\\
= \Big( \sum_{i=1}^N \ind{(i_{j-1},i_{j}]}(x) \frac{\int_{(i_{j-1},i_{j}]} u(y)dy}{i_{j}-i_{j-1}}  \Big) +
\sum_{i=1}^N \ind{(i_{j-1},i_{j}]}(x) \Big( u(x) - \frac{\int_{(i_{j-1},i_{j}]} u(y)dy}{i_{j}-i_{j-1}}  \Big) .
\end{multline}
We prove the following theorem.
\begin{thm}
Let $W \in \mathcal{P}_{\vi\times\vi} $. There exists a unique solution $u \in C([0,T],L^2([0,1]))$ 
of \eqref{eq:GraphonChaucy}
for any initial condition $u_0\in L^2([0,1])$ and $T>0$.
\end{thm}
\begin{proof} \textbf{Uniqueness}. 
Multiplying both terms of \eqref{eq:GraphonChaucy} by $u(x,t)$ and integrating on $x$ we obtain

\[
\frac{d}{dt} \frac{\int_0^1 u^2(x,t) dx}{2}
= \int_0^1 \Big( \int_{[0,1]} W(x,y) (u(y,t) - u(x,t)) dy \Big) u(x,t) dx
\leq 2 C {\int_0^1 u^2(x,t) dx}
\]
which implies the uniqueness by Gronwall Lemma, since $\|u_t\|_2^2 \leq \|u_0\|_2^2 \exp(2Ct)$.

\smallskip

\noindent\textbf{Existence}.
The superposition principle allows us to find the solutions of the Graphon Chauchy problem \eqref{eq:GraphonChaucy}
on each component of the orthogonal decomposition 
$L^2(0,1) = \mathcal{P}_{\vi} \oplus_{j=1}^N L^2_0 ((i_{j-1},i_{j}]) $
given above. In fact, it is sufficient to use
\eqref{eq:superposit} to decompose the initial condition, if $u_0\in L^2([0,1])$. 

\medskip

\emph{Solution starting in $\mathcal{P}_\vi$.}
When $u_0\in \mathcal{P}_\vi$, it may be written as $u_0 (x)= \sum_{j=1}^N u_{0j} \ind{( i_{j_1-1},i_{j_1}]} (x)$.
Denote by $\vu_0 = (u_{01},\ldots,u_{0N})^\top$ and by $\ind{\vi} = (\ind{(i_{0},i_{1}]}, \ind{(i_{1},i_{2}]}, \ldots ,\ind{(i_{N-1},i_{N}]})^\top$,
so that $u_0(x) = \ind{\vi}(x)^\top \vu_0$.
We claim that the (unique) solution of \eqref{eq:GraphonChaucy} belongs to 
$\mathcal{P}_\vi$ for any $t \geq 0$ and it has the form

$$
u(x,t) = \ind{\vi}(x)^\top \vu(t) = \ind{\vi}(x)^\top \exp(-\Delta_{W_\vi} t) \vu_0,
$$ where
$\Delta_{W_\vi}$ is the Laplacian of the $N\times N$ matrix $W_\vi$ given by
\begin{equation*}
[W_\vi]_{j_1j_2}
= b_{j_1j_2}({i_{j_2}-i_{j_2-1}})
= \int_{[0,1]} \ind{((i_{j_2-1},i_{j_2}]}(y) W(x,y) dy  \qquad \text{for any }x \in \ind{((i_{j_1-1},i_{j_1}]} 
.
\end{equation*}
Note that, in matrix form

\begin{equation}\label{eq:Wast}
\ind{\vi}(x)^\top diag(W_\vi\vone) = \Big( \int_{[0,1]} W(x,y) dy \Big) \ind{\vi}(x)^\top ,
\qquad
\ind{\vi}(x)^\top W_\vi = \int_{[0,1]} W(x,y) \ind{\vi}(y)^\top  dy.
\end{equation}
In fact:
\begin{itemize}
\item by definition, $u(\cdot,t) \in \mathcal{P}_\vi$ for any $t \geq 0$ and $u(x,0) = u_0(x)$;
\item by definition of exponential matrix and by \eqref{eq:Wast},
the equation \eqref{eq:GraphonChaucy} is satisfied, since 

\[
\begin{aligned}
\frac{ \partial  u}{\partial t}(x,t)
& =
\ind{\vi}(x)^\top  (- \Delta_{W_\vi} ) \exp(-\Delta_{W_\vi} t) \vu_0
=
\ind{\vi}(x)^\top ( {W_\vi} - diag({W_\vi}\vone) ) \exp(-\Delta_{W_\vi} t)  \vu_0
\\
& =
\int_{[0,1]} W(x,y) \Big[\ind{\vi}(y)^\top \exp(-\Delta_{W_\vi} t) \vu_0\Big] dy
-
\int_{[0,1]} W(x,y) \Big[\ind{\vi}(x)^\top \exp(-\Delta_{W_\vi} t) \vu_0\Big] dy
\\
& =
\int_{[0,1]} W(x,y) (u(y,t) - u(x,t)) dy
.
\end{aligned}
\]
\end{itemize}
\begin{rem}
Note that the Laplacian of a square matrix does not depend on the values on its diagonal.
As an expected consequence, $(W_{\vi})_{j j}$ does not affect the action of \eqref{eq:GraphonChaucy}
on $\mathcal{P}_{\vi}$, but it has an effect on $\oplus_{j=1}^N L^2_0 (i_{j-1},i_{j})$.
\end{rem}

\medskip

\emph{Solution starting in $L^2_0 ((i_{j-1},i_{j}])$}. 
Fixed $j \in \{1,\ldots,N\}$, for $u_0\in L^2_0 ((i_{j-1},i_{j}])$, we claim that the (unique) solution of \eqref{eq:GraphonChaucy} has the form
$ u(x,t) = u_0(x) \exp(-\mu t) $, where $\mu = \sum_{j_2=1}^N b_{j\,j_2} (i_{j_2}-i_{j_2-1}) = (W_\vi\vone)_j $.
In fact:
\begin{itemize}
\item since $u_0 \in L^2_0 ((i_{j-1},i_{j}])$ and $u_0(\cdot) = u_0 \ind{(i_{j-1},i_{j}]}(\cdot) = 
u(x,0)$,
then $u(\cdot,t) \in L^2_0 ((i_{j-1},i_{j}])$ for any $t \geq 0$;
\item the equation \eqref{eq:GraphonChaucy} is satisfied, again since $u_0(\cdot) = u_0(\cdot) \ind{(i_{j-1},i_{j}]}(\cdot)$,
and hence

\begin{align*}
\frac{ \partial  u}{\partial t}(x,t) & =
- \mu u_0(x) \exp(-\mu t) \\
& = \exp(-\mu t) \Big(
0 - u_0(x) \sum_{j_2=1}^N b_{j\,j_2} (i_{j_2}-i_{j_2-1})
\Big)
\\
& =
\exp(-\mu t) \Big(
\Big[
\sum_{j_1=1}^N b_{j_1j} \ind{( i_{j_1-1},i_{j_1}]}(x)
\Big]
\int_{i_{j-1}}^{i_{j}} u_0(y)dy
\\
& \qquad
- u_0(x) \int_{[0,1]}
\Big[
\ind{(i_{j-1},i_{j})}(x) \sum_{j_2=1}^N b_{j\,j_2} \ind{(i_{j_2-1},i_{j_2}]} (y)
\Big]
dy\Big)
\\
& =
\exp(-\mu t) \Big(
\int_0^1 W(x,y) u_0(y)\ind{(i_{j-1},i_{j}]}(y) dy - u_0(x) \ind{(i_{j-1},i_{j}]}(x) \int_0^1 W(x,y) dy
\Big)
\\
& =
\exp(-\mu t) \Big(
\int_0^1 W(x,y) (u_0(y)-u_0(x)) dy
\Big)
\,=\,
\int_{[0,1]} W(x,y) (u(y,t) - u(x,t)) dy.\qedhere
\end{align*} 
\end{itemize}
\end{proof}

\section{Dynamics for symmetric cases with $N=3$ groups}
Let $N = 3$ and $\vi = (i_0=0, i_1,i_2,i_3=1)$ and assume
$W_{\vi}$ to be symmetric. Then

\[
{W_{\vi}} = \begin{pmatrix}
 a_{11} & a_{12} & a_{13}
\\
a_{12} & a_{22} & a_{23} \\
a_{13} & a_{23} & a_{33}
\end{pmatrix},
\qquad
\Delta_{W_{\vi}} = \begin{pmatrix}
 a_{12}+a_{13} & -a_{12} & -a_{13}
\\
-a_{12} & a_{23}+a_{12} & -a_{23} \\
-a_{13} & -a_{23} & a_{13}+a_{23}
\end{pmatrix}
\]
where $a_{j_1j_2}= b_{j_1j_2} (i_{j_2}-i_{j_2-1})$ and $0<j_1\leq j_2\leq3$.
Denoting by

\[
\Delta = \sqrt{a_{12}^2+a_{13}^2+a_{23}^2 -a_{12}a_{13} - a_{12}a_{23} - a_{13}a_{23}}
= \sqrt{\frac{(a_{12}-a_{13})^2+(a_{12}-a_{23})^2 +(a_{13} - a_{23})^2}{2}}
\]
we get the eigenvalues of $\Delta_{W_{\vi}}$:

\[
\lambda_1 = 0,
\quad
\lambda_{2} = a_{12}+a_{13}+a_{23} + \Delta,
\quad
\lambda_{3} = a_{12}+a_{13}+a_{23} - \Delta
\]
and the three corresponding orthogonal eigenvectors

\[
\vv_1 =
\begin{pmatrix}
1
\\
1
\\
1
\end{pmatrix},
\vv_2 =
\begin{pmatrix}
a_{23} - a_{12} - \Delta
\\
a_{12} - a_{13} + \Delta
\\
a_{13} - a_{23}
\end{pmatrix},
\vv_3 =
\begin{pmatrix}
a_{23} - a_{12} + \Delta
\\
a_{12} - a_{13} - \Delta
\\
a_{13} - a_{23}
\end{pmatrix},
\qquad
V = (\vv_1 | \vv2 | \vv_3)
\]
so that $\Delta_{W_{\vi}} V = V diag(\lambda_1,\lambda_2,\lambda_3) $.
Moreover, any symmetric matrix has orthogonal
eigenspaces, and hence
\(
V^\top V =
diag(
\|\vv_1\|^2,
\|\vv_2\|^2,
\|\vv_3\|^2
)
\).
We note that $\exp( -\Delta_{W_{\vi}} t) = V
diag(1,e^{-\lambda_2t},e^{-\lambda_2t}) V^{-1}$, where $V^{-1} =
diag(
1/\|\vv_1\|^2,
1/\|\vv_2\|^2,
1/\|\vv_3\|^2
)
 V^{\top}$.

Now, let $u_0(x) \in L^2([0,1])$ be the initial condition; we can immediately find the solution
$u(x,t) = \ind{\vi}(x)^\top \vu(t) $
of the component in $\mathcal{P}_{\vi}$ with

\begin{equation}\label{eq:parties}
\begin{aligned}
\vu_0 = &
\begin{pmatrix}
\frac{\int_{i_0=0}^{i_1} u_0(y) dy}{i_{1}-i_{0}}
&
\frac{\int_{i_1}^{i_2} u_0(y) dy}{i_{2}-i_{1}}
&
\frac{\int_{i_2}^{i_3=1} u_0(y) dy}{i_{3}-i_{2}}
\end{pmatrix}^\top,
\\
\vu (t)
& = \exp( -\Delta_{W_{\vi}} t) \vu_0
\\
&
=
V
diag(1,e^{-\lambda_2t},e^{-\lambda_2t}) diag(
1/\|\vv_1\|^2,
1/\|\vv_2\|^2,
1/\|\vv_3\|^2
)
 V^{\top} \vu_0.
\end{aligned}
\end{equation}
The vector $\vu_0$ gives that barycenter of the three groups' opinion.
Note that the ``barycenter'' $b(t)$ of the whole system does not change with $t$, since $\vv_1=\vone^\top$ is a (left)
kernel vector of $\Delta_{W_{\vi}}$, and hence of $\exp( -\Delta_{W_{\vi}} t)$. In fact, let
$b(t) = \frac{\vone^\top}{3} \vu(t)$. Since $\frac{\vone^\top}{3} ^\top V = (1 ,0, 0)^\top$,
we have $b(t) = \frac{1}{\|\vv_1\|^2} \vv_1^\top \vu_0 = \frac{\vone^\top}{3} \vu_0 = b(0)$.

By superposing the solution with that of $L^2_0 ((i_{j-1},i_{j}])$, we find the final solution:

\begin{align*}
&\text{group }1: && \forall x \in [0,t_1], &&
u(x,t) = u_1 (t) + ( u_0(x) - u_{01}) \exp( - ( a_{11} + a_{12} + a_{13}) t )
\\
&\text{group }2: && \forall x \in (t_1,t_2], &&
u(x,t) = u_2 (t) + ( u_0(x) - u_{02}) \exp( - ( a_{12} + a_{22} + a_{23}) t )
\\
&\text{group }3: && \forall x \in (t_2,1], &&
u(x,t) = u_3 (t) + ( u_0(x) - u_{03}) \exp( - ( a_{13} + a_{23} + a_{33}) t ).
\end{align*}
\subsection{Right, center, and left parties}
Assume $0 \leq \min(a_{12},a_{23})$ and $a_{12}+a_{23}>0$. In other words,
we think that the second group acts as a political``center'', or mediator, that tries to attract the other groups toward
a common center. The other groups may have a conflict relationship ($a_{13}<0$) and the strength
of this conflict may lead to different scenarios. 
We start by noticing that the assumption on the second group implies that, when  $a_{13}<0$,
\(
\Delta^2 \geq \frac{ a_{13}^2 + (a_{12}-a_{23})^2 + a_{13}^2}{2}
\), and then in this case $\Delta > |a_{13}| $, which implies that $\lambda_2 > 0$ (that is trivially true for $a_{13}\geq 0$).
Now, $\lambda_2 > 0$ forces the initial component in $\vv_2 $
to collapse to $0$ (see \eqref{eq:parties}). The possible different scenarios depend on the sign of $\lambda_3$, which is
always lower than $\lambda_2$.\\
{\it Dynamics of the barycenter of the three groups as function of $\lambda_3$.}
We have, for what regards the dynamics on $\mathcal{P}_{\vi}$, with respect to $a_{13}$
\begin{itemize}
\item[$\lambda_3>0$:] this result is achieved when $a_{13} > -\tfrac{a_{12}a_{23}}{a_{12}+a_{23}}$.
With $\lambda_3>0$,  also the component in
$\vv_3$ is forced to tend to $0$ by
\eqref{eq:parties}. Then the limit invariant subspace
of $\mathcal{P}_{\vi}$ is the constant one: $\vu(t) \to b(0)\vone$, i.e.\ all the three barycenter collapse into a single
point determined by the initial conditions;
\item[$\lambda_3 = 0$:] when $a_{13} = -\tfrac{a_{12}a_{23}}{a_{12}+a_{23}}$, by substituting we get
$\lambda_2= 2 (a_{12}+a_{13}+a_{23}) = 2\tfrac{a_{12}^2+a_{23}a_{12}+a_{23}^2}{a_{12}+a_{23}}$, and 

\[
\vv_1 =
\begin{pmatrix}
1
\\
1
\\
1
\end{pmatrix},
\vv_2 =
\begin{pmatrix}
-2a_{12} + \tfrac{a_{12}a_{23}}{a_{12}+a_{23}}
\\
2 a_{12}+ a_{23}
\\
- a_{23} -  \tfrac{a_{12}a_{23}}{a_{12}+a_{23}}
\end{pmatrix},
\vv_3 =
\begin{pmatrix}
2 a_{23} + \tfrac{a_{12}a_{23}}{a_{12}+a_{23}}
\\
- a_{23} -2 \tfrac{a_{12}a_{23}}{a_{12}+a_{23}}
\\
- a_{23} + \tfrac{a_{12}a_{23}}{a_{12}+a_{23}}
\end{pmatrix}.
\]
The fact that $\lambda_3=0$ leaves the component on $\vv_3$ unchanged as that on $\vv_1$, and then by \eqref{eq:parties} there will be the formation of three different barycenter of the groups in

\[
\vu_{\infty} = \lim_{t \to \infty} \vu(t) = \vone\frac{\vone^\top \vu_0 }{3} +\vv_3 \frac{\vv_3^\top \vu_0 }{\|\vv_3\|^2};
\]
\item[$\lambda_3 < 0$]: for $a_{13} < -\tfrac{a_{12}a_{23}}{a_{12}+a_{23}}$,
the barycenter of the whole system remains constant (the component on $\vv_1$), while the barycenter of the three clusters diverge
by \eqref{eq:parties}.
\end{itemize}
{\it Dynamics inside each group with respect to its own barycenter.}\\
The dynamic in each $L^2_0 ((i_{j-1},i_{j}])$ depends on the sign of $\mu_j = (W_\vi\vone)_j = \sum_{j_2=1}^3 a_{j j_2}$. In fact,
as a consequence of the proof above,
\(
( u(x,t) - u_j(t)) \ind{(i_{j-1},i_{j}]} (x) = \exp( - \mu_j t ).
\)
As an example, for the second group (the other are similar), we have $\mu_2 = a_{12} + a_{22} + a_{23} $. Hence
\begin{description}
\item[$\mu_2 >0$:] the second group will collapse on its barycenter;
\item[$\mu_2 =0$:] in this equilibrium phase, the points of the second group
are translated rigidly with their barycenter. This happens when the repulsion inside the center party equalizes the action of this party on the opposite ones: $a_{22} = - (a_{12}+a_{23} )$;
\item[$\mu_2 <0$:] the repulsion inside the center party is so big ($a_{22} < - (a_{12}+a_{23} )$) that its points will explode exponentially fast around its barycenter.
\end{description}


\begin{thebibliography}{10}
\expandafter\ifx\csname url\endcsname\relax
  \def\url#1{\texttt{#1}}\fi
\expandafter\ifx\csname urlprefix\endcsname\relax\def\urlprefix{URL }\fi
\expandafter\ifx\csname href\endcsname\relax
  \def\href#1#2{#2} \def\path#1{#1}\fi

\bibitem{IJMOC2020}
H.~Noorazar, K.~Vixie, A.~Talebanpour, Y.~Hu, From classical to modern opinion
  dynamics, International Journal of Modern Physics C 31~(07) (2020) 2050101.
\newblock \href {http://dx.doi.org/10.1142/S0129183120501016}
  {\path{doi:10.1142/S0129183120501016}}.

\bibitem{PT2013}
L.~Pareschi, G.~Toscani,
  \href{https://books.google.it/books?id=lct7AgAAQBAJ}{Interacting Multiagent
  Systems: Kinetic equations and Monte Carlo methods}, OUP Oxford, 2013.
\newline\urlprefix\url{https://books.google.it/books?id=lct7AgAAQBAJ}

\bibitem{sznajd2000opinion}
K.~Sznajd-Weron, J.~Sznajd, Opinion evolution in closed community,
  International Journal of Modern Physics C 11~(06) (2000) 1157--1165.

\bibitem{clifford1973model}
P.~Clifford, A.~Sudbury, A model for spatial conflict, Biometrika 60~(3) (1973)
  581--588.

\bibitem{galam1986}
S.~Galam, Majority rule, hierarchical structures, and democratic
  totalitarianism: A statistical approach, Journal of Mathematical Psychology
  30~(4) (1986) 426--434.

\bibitem{Degroot1974118}
M.~Degroot, Reaching a consensus, Journal of the American Statistical
  Association 69~(345) (1974) 118 – 121.
\newblock \href {http://dx.doi.org/10.1080/01621459.1974.10480137}
  {\path{doi:10.1080/01621459.1974.10480137}}.

\bibitem{Friedkin1990193}
N.~Friedkin, E.~Johnsen, Social influence and opinions, The Journal of
  Mathematical Sociology 15~(3-4) (1990) 193 – 206.
\newblock \href {http://dx.doi.org/10.1080/0022250X.1990.9990069}
  {\path{doi:10.1080/0022250X.1990.9990069}}.

\bibitem{Hegselmann2002}
R.~Hegselmann, U.~Krause, Opinion dynamics and bounded confidence: Models,
  analysis and simulation, JASSS 5~(3).

\bibitem{Deffuant2000}
G.~Deffuant, D.~Neau, F.~Amblard, G.~Weisbuch, Mixing beliefs among interacting
  agents, Advances in Complex Systems 3~(1/4) (2000) 87--98.
\newblock \href {http://dx.doi.org/10.1142/S0219525900000078}
  {\path{doi:10.1142/S0219525900000078}}.

\bibitem{MYE2019}
M.~Ye, Opinion Dynamics and the Evolution of Social Power in Social Networks,
  Springer Cham, Switzerland, 2019.

\bibitem{Noorazar2020}
H.~Noorazar, Recent advances in opinion propagation dynamics: a 2020 survey,
  European Physical Journal Plus 135~(6).
\newblock \href {http://dx.doi.org/10.1140/epjp/s13360-020-00541-2}
  {\path{doi:10.1140/epjp/s13360-020-00541-2}}.

\bibitem{Benf6}
A.~Benfenati, V.~Coscia, Modeling opinion formation in the kinetic theory of
  active particles {I}: spontaneous trend, Ann. Univ. Ferrara Sez. VII Sci.
  Mat. 60~(1) (2014) 35--53.
\newblock \href {http://dx.doi.org/10.1007/s11565-014-0207-2}
  {\path{doi:10.1007/s11565-014-0207-2}}.

\bibitem{Benf7}
A.~Benfenati, V.~Coscia, Nonlinear microscale interactions in the kinetic
  theory of active particles, Appl. Math. Lett. 26~(10) (2013) 979--983.
\newblock \href {http://dx.doi.org/10.1016/j.aml.2013.04.007}
  {\path{doi:10.1016/j.aml.2013.04.007}}.

\bibitem{ANT07}
G.~Aletti, G.~Naldi, G.~Toscani,
  \href{https://doi.org/10.1137/060658679}{First‐order continuous models of
  opinion formation}, SIAM Journal on Applied Mathematics 67~(3) (2007)
  837--853.
\newblock \href {http://arxiv.org/abs/https://doi.org/10.1137/060658679}
  {\path{arXiv:https://doi.org/10.1137/060658679}}, \href
  {http://dx.doi.org/10.1137/060658679} {\path{doi:10.1137/060658679}}.
\newline\urlprefix\url{https://doi.org/10.1137/060658679}

\bibitem{Borgs2008}
C.~Borgs, J.~Chayes, L.~Lov\'{a}sz, V.~S\'{o}s, K.~Vesztergombi, Convergent
  sequences of dense graphs i: Subgraph frequencies, metric properties and
  testing, Advances in Mathematics 219~(6) (2008) 1801 – 1851.
\newblock \href {http://dx.doi.org/10.1016/j.aim.2008.07.008}
  {\path{doi:10.1016/j.aim.2008.07.008}}.

\bibitem{Borgs2012}
C.~Borgs, J.~Chayes, L.~Lov\'{a}sz, V.~S\'{o}s, K.~Vesztergombi, Convergent
  sequences of dense graphs ii. multiway cuts and statistical physics, Annals
  of Mathematics 176~(1) (2012) 151 – 219.
\newblock \href {http://dx.doi.org/10.4007/annals.2012.176.1.2}
  {\path{doi:10.4007/annals.2012.176.1.2}}.

\bibitem{Gao20}
S.~Gao, P.~Caines,
  \href{https://www.scopus.com/inward/record.uri?eid=2-s2.0-85080134638&doi=10.1109%2fTAC.2019.2955976&partnerID=40&md5=676b787a360781ac1d1851b32d7d6548}{Graphon
  control of large-scale networks of linear systems}, IEEE Transactions on
  Automatic Control 65~(10) (2020) 4090--4105, cited By 15.
\newblock \href {http://dx.doi.org/10.1109/TAC.2019.2955976}
  {\path{doi:10.1109/TAC.2019.2955976}}.
\newline\urlprefix\url{https://www.scopus.com/inward/record.uri?eid=2-s2.0-85080134638&doi=10.1109%2fTAC.2019.2955976&partnerID=40&md5=676b787a360781ac1d1851b32d7d6548}

\bibitem{Benf9}
G.~Aletti, A.~Benfenati, G.~Naldi, A semi-supervised reduced-space method for
  hyperspectral imaging segmentation, Journal of Imaging 7~(12).
\newblock \href {http://dx.doi.org/10.3390/jimaging7120267}
  {\path{doi:10.3390/jimaging7120267}}.

\bibitem{Benf11}
G.~Aletti, A.~Benfenati, G.~Naldi, A semiautomatic multi-label color image
  segmentation coupling {D}irichlet problem and colour distances, Journal of
  Imaging 7~(10).
\newblock \href {http://dx.doi.org/10.3390/jimaging7100208}
  {\path{doi:10.3390/jimaging7100208}}.

\bibitem{Benf10}
G.~Aletti, A.~Benfenati, G.~Naldi, Graph, spectra, control and epidemics: An
  example with a {SEIR} model, Mathematics 9~(22).
\newblock \href {http://dx.doi.org/10.3390/math9222987}
  {\path{doi:10.3390/math9222987}}.

\bibitem{KVM2018}
D.~Kaliuzhnyi-Verbovetskyi, G.~Medvedev, The mean field equation for the
  {K}uramoto model on graph sequences with non-lipschitz limit, SIAM Journal on
  Mathematical Analysis 50~(3) (2018) 2441--2465.
\newblock \href {http://dx.doi.org/10.1137/17M1134007}
  {\path{doi:10.1137/17M1134007}}.

\end{thebibliography}
\end{document}